\documentclass[]{imsart}

\RequirePackage{amsthm,amsmath}
\RequirePackage[numbers]{natbib}
\RequirePackage[colorlinks,citecolor=blue,urlcolor=blue]{hyperref}

\startlocaldefs
\numberwithin{equation}{section}
\theoremstyle{plain}
\newtheorem{thm}{Theorem}[section]
\newtheorem{proposition}{Proposition}[section]
\newtheorem{remark}{Remark}[section]
\endlocaldefs

\usepackage{natbib}
 \bibpunct[, ]{(}{)}{,}{a}{}{,}%
 \def\bibsep{\smallskipamount}%
\usepackage{amsmath}
\usepackage{amssymb}
\usepackage{upgreek}
\usepackage{bbold} 

\usepackage{graphicx}
\usepackage{tikz}
\usetikzlibrary{shapes,backgrounds,patterns,arrows}
\usepackage{pgfplots}
\pgfplotsset{compat=1.12}
\usepgfplotslibrary{groupplots}
\usepackage{subfigure}
\usepackage{appendix}

\def\barmu{\bar\mu}

\def\({{\Bigl(}}
\def\){{\Bigr)}}

\newcommand{\vertexUnexplored}{\tikz[baseline=-0.55ex]\draw[black,fill=white,radius=0.4ex] (0,0) circle;}
\newcommand{\vertexActive}{\tikz[baseline=-0.55ex]\draw[black,fill=black,radius=0.4ex] (0,0) circle;}
\newcommand{\vertexBlocked}{\tikz[baseline=-0.55ex]\draw[black,fill=red,radius=0.4ex] (0,0) circle;}

\usepackage{glossaries} 

\newglossaryentry{Erdos}{
name=\ensuremath{\mbox{ER}_n(\lambda/n)},
description={Erdos-R\'enyi graph}
}

\newglossaryentry{CM}{
name=\ensuremath{\mbox{CM}_n(\bar d^{(n)})},
description={Configuration Model graph}
}

\begin{document}

\begin{frontmatter}
\runtitle{Degree-Greedy Algorithms on Large Random Graphs}

\title{Sequential Algorithms and Independent Sets Discovering on Large Sparse Random Graphs}


%
%
%
%
%
%
%

\author{Paola Bermolen$^1$,Matthieu Jonckheere$^2$, Federico Larroca$^1$, Manuel Saenz$^2$}\\

\affiliation{}{$^1$ Facultad de Ingeniería, Universidad de la República, Uruguay}
\affiliation{}{$^2$Instituto de C\'alculo, UBA/CONICET, Buenos Aires, Argentina}

\address{paola@fing.edu.uy, matthieu.jonckheere@gmail.com, flarroca@fing.edu.uy, manuel.saenz@gmail.com}
\thanks{This work was partially support by Stic Amusd GENE project}

\begin{abstract}
 Computing the size of maximum independent sets is a NP-hard problem for fixed graphs. Characterizing and designing efficient algorithms to estimate this independence number for random graphs are notoriously difficult and still largely open issues. In a companion paper, we showed that a low complexity degree-greedy exploration is actually asymptotically optimal on a large class of sparse random graphs. Encouraged by this result, we present and study two variants of sequential exploration algorithms: static and dynamic degree-aware explorations. We derive hydrodynamic limits for both of them, which in turn allow us to compute the size of the resulting independent set. Whereas the former is simpler to compute, the latter may be used to arbitrarily approximate the degree-greedy algorithm. Both can be implemented in a distributed manner. The corresponding hydrodynamic limits constitute an efficient method to compute or bound the independence number for a large class of sparse random graphs. As an application, we then show how our method may be used to estimate the capacity of a large 802.11-based wireless network. We finally consider further indicators such as the fairness of the resulting configuration, and show how an unexpected trade-off between fairness and capacity can be achieved. 
\end{abstract}


\end{frontmatter}




%


\section{Introduction}

Given a graph, an \emph{independent set} is a subset of its vertices where no pair of them are connected to each other. An independent set is said to be \emph{maximal} if no more vertices may be added to it, whereas it is said to be \emph{maximum} if no other independent set is larger in size (and this size is usually referred to as the \emph{independence number} of the graph).

Computing the independence number of random graphs is a theoretically important and difficult problem, which finds applications in several areas, e.g. in biology, physics and engineering. For instance, by considering a wireless network where nodes (or links) are represented as vertices in the graph, and two nodes (or links) are connected by an edge when they cannot transmit simultaneously as a consequence of the medium access mechanism (i.e., the so-called interference graph), the authors of \cite{laufer2016capacity} proved that the capacity of a 802.11-based wireless network is given by the maximum size of its independent sets (see also~\cite{soung2010back}). 

However, computing the independence number of a given graph is well-known to be an NP-hard problem. Our main focus here is thus finding efficient methods to compute (up to vanishing errors) the independence number of large sparse random graphs. Building on recent results~\cite{bermolen2017jamming,janson2017greedy} we study the maximal independent sets obtained by sequential algorithms in sparse random graphs. For example, one such algorithm could be to successively select uniformly at random vertices of the graph one at a time, and add them to the independent set if none of their neighbors already belong to it. This exploration procedure is known as the \emph{Greedy Algorithm}. 

Clearly, the size of the independent set obtained by the greedy algorithm is in general smaller than the independence number. Larger independent sets may be obtained by considering the degree of the vertices chosen at each step. In this sense, a classical algorithm to approximate the independence number is the \emph{Degree-Greedy Algorithm}; i.e.\ an exploration giving at each step full priority to those vertices that have the smallest degree towards those that were not explored yet.

We first review in detail the characterization of the size of the independent set obtained by sequential random explorations on sparse random graphs (Erd{\"o}s-R\'enyi and Configuration Models), underlining that few results are known on the independence number of these graphs.
We give a particular attention to recently obtained functional laws of large numbers when the size of the graph tends to infinity. 
We recently proved in a companion paper \cite{jonckheere2018asymptotic} that the degree-greedy algorithm is asymptotically optimal for a large class of sparse random graphs.  We discuss here various examples of such graphs (as well as examples of graphs where degree-greedy is not asymptotically optimal).

Our contribution is then four-fold. 
Our first contribution consists in showing how one can deploy and evaluate a decentralized version of the degree-greedy algorithm,
by generalizing the hydrodynamic limits known for random sequential algorithms~\cite{bermolen2017jamming,janson2017greedy} to a large class of degree-aware sequential algorithms. More specifically, we define two types of mechanisms: one that takes into account only the initial degree (which we will be calling \emph{Static Degree-aware Exploration Algorithms}), and one which takes into account the degree adaptively (\emph{Dynamic Degree-aware Exploration Algorithms}).

Regarding the static degree-aware exploration algorithm, although sub-optimal, the resulting hydrodynamic limit is solvable exactly, i.e., in closed form. This allows us, for instance, to calculate the size of the corresponding independent set but also to calculate the limit of the proportion of vertices with a given degree in the independent set.
On the other hand, the dynamic degree-aware exploration algorithms may be used to approximate to an arbitrary precision the performance of the degree-greedy algorithm. Although the resulting hydrodynamic limit is not solvable analytically and we have to resort to numerical estimations, this constitutes an efficient and close estimation of the independence number for a large family of random graphs. As an example of application, we then compare the capacity of 802.11-based networks and the proposed estimation. By means of several simulations, we show that the performance of the latter is very close to the optimal one, even when the underlying random graph does not belong to this family.  It is important to note that while the degree-greedy cannot be implemented in a distributed manner, dynamic and static degree-aware explorations can, meaning that this schemes could in principle be implemented as a communication protocol. 

A second contribution consists in studying numerically the benefits of an algorithm that combines both strategies (degree-greedy and Glauber) in order to achieve quasi-optimal results much faster in all cases (high and low connectivity).

As a third contribution, we show that one of the modified sequential algorithms allows to reach interesting tradeoffs between fairness (equality of chances to access the communication channel for nodes of different degrees) and efficiency: fairness can be significantly improved while losing only little capacity.   

Our last contribution consists in comparing numerically our findings for large sparse random graphs with all pre-existing results which have not yet been compared between them, even in the simplest Erd{\"o}s R\'enyi case. Besides, we show how these insights allow to predict performance indicators of real large wireless networks better, for instance, than stochastic geometry.

The rest of the article is structured as follows. In Section ~\ref{sec:state_of_the_art} we define our model and present previous results on sequential algorithms for the exploration of random graphs. In particular, we summarize known results for the greedy and degree-greedy algorithms, as well as for the characterization of the maximum independent set. In Section ~\ref{sec:optimality} we recall conditions for the optimality of the degree-greedy algorithm, obtained in a companion paper. In Section ~\ref{sec:generalization_seq_algorithms} we introduce and analyze two variants of sequential algorithms. There, we prove hydrodynamic limits and characterize the associated independent set sizes. In Section ~\ref{sec:simulaciones} we conduct a performance analysis through several simulations on different random graphs. Finally, in Section ~\ref{sec:maxywifi} we address the problem of approximating  the capacity of large wireless network. In this section, we also study the trade off between fairness and efficiency.\\

A shorter version of this work was published at the proceedings of the \emph{36th International Symposium on Computer Performance, Modeling, Measurements and Evaluation 2018} (\cite{PER}). The present article has several important extensions of our previous work:

\begin{itemize}
    \item it presents a more detailed bibliographical review of the problem of finding maximal independent sets on random graphs,
    \item it describes in depth the relationship of our work with previous results,
    \item the proof of key results are presented in more detail (and the connection to previous results explained),
    \item a discussion on the possible combination of sequential and dynamical (CSMA-like) algorithms and simulations showing their improvements in terms of convergence times are included,
    \item it incorporates several simulations providing interesting insights (for instance, comparing previous results found in the literature, but also considering real-life graphs),
    \item and it also counts with an analysis of the impact of different strategies on the fairness measures of the network.
\end{itemize}


\section{Model and previous works}\label{sec:state_of_the_art}

In the sequel, $G^{(n)} = (V,E)$ denotes a (possibly random) \emph{graph} that consists of a set of $n$ vertices $V = \{ 1, \ldots, n \}$, and a set of undirected edges $E \subseteq V \times V$. 
 
As mentioned in the previous section, given a graph $G^{(n)}$, an \emph{independent set} is a subset of vertices $A \subseteq V$ where for every pair $v, w \in A$ we have that $\{ v, w \} \notin E$ (i.e.\ no pair of vertices are connected to each other). An independent set is said to be \emph{maximal} if it is not the subset of a larger independent set; and \emph{maximum} if there is no other independent set of larger size. The size of a maximum independent set is usually referred as the \emph{independence number} of a graph and denoted by $\alpha(G^{(n)})$.
 
Two types of explorations have been considered both in the random graph literature and for the modelling of wireless networks: \emph{sequential} and \emph{dynamical algorithms}. First, we will discuss the sequential case, where nodes are successively added to an independent set until a maximal one is reached. We will present some results on two sequential exploration examples: greedy and degree-greedy algorithms. Afterwards, we will treat the dynamical case where active nodes can reverse their state to unexplored after some random time. More specifically, we will review some aspects of the so-called Glauber dynamics~\cite{vigoda2001note}.


\subsection{Sequential exploration}

\begin{figure}[ht]
\centering
\subfigure[$k=1$]{
    \includegraphics[scale=1.2]{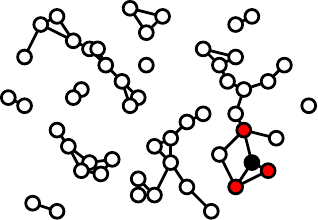}
    \label{fig:subfig1}
}
\subfigure[$k=3$]{
    \includegraphics[scale=1.2]{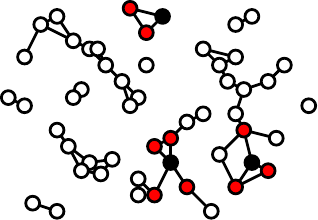} 
    \label{fig:subfig2}
}
\subfigure[$k=k_n^*$]{
    \includegraphics[scale=1.2]{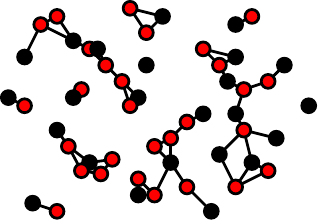}
    \label{fig:subfig3}
}
\caption{An exploration process on a fixed graph, that stopped at some time $k_n^*$. Image courtesy of Jaron and Sanders.}
\label{fig:Exploration_process}
\end{figure}

In what follows, random sequential algorithms will be used to find maximal independent sets and analyze them. At any step $k = 0, 1, 2, \ldots$  of the algorithms discussed, we will consider that each vertex is either \emph{active} \vertexActive, \emph{blocked} \vertexBlocked, or \emph{unexplored} \vertexUnexplored\ (see Fig.\ \ref{fig:Exploration_process}). Accordingly, the set of vertices will be split into three components: the set of active vertices $\mathcal{A}_k$, the set of blocked vertices $\mathcal{B}_k$, and the set of unexplored vertices $\mathcal{U}_k$. Active vertices may be taken to correspond to nodes that transmit, blocked vertices to nodes that cannot transmit because they are impeded of doing so by a neighboring active node, and unexplored vertices to nodes that are in neither of these two states. At any step $k$, the active vertices will be the ones that belong to the independent set constructed by the algorithm.

A typical sequential exploration algorithm works in the following way. Initially, it sets $\mathcal{U}_0 = V$, $\mathcal{A}_0 = \emptyset$ and $\mathcal{B}_0 = \emptyset$. To explore the graph, at the $k+1$-th step it selects a vertex $v_{k+1}\in\mathcal{U}_k$ (possibly taking into account its current or past degree towards other unexplored vertices), and changes its state into active. After this, it takes all of its unexplored neighbors, i.e.\ the set $\mathcal{N}_{v_{k+1}} = \{ w \in \mathcal{U}_k | (v_{k+1},w) \in E \}$, and changes their states into blocked. This means that the resulting set of vertices will be given by $\mathcal{U}_{k+1} = \mathcal{U}_k \backslash \{ v_{k+1} \cup \mathcal{N}_{v_{k+1}} \}$, $\mathcal{A}_{k+1} = \mathcal{A}_{k+1} \cup \{ v_{k+1} \}$ and $\mathcal{B}_{k+1} = \mathcal{B}_k \cup \mathcal{N}_{v_{k+1}}$. As mentioned before, at each step the set of active vertices defines an independent set. The algorithm keeps repeating this procedure until the step $k_n^*$ in which all vertices are either active or blocked (or equivalently $\mathcal{U}_{k_n^*} = \emptyset$). An example of an exploration process on a fixed graph is depicted in Fig. \ref{fig:Exploration_process}. 

The set of active vertices at step $k_n^*$ then defines a maximal independent set. The proportion of vertices contained in this independent set $|\mathcal{A}_{k_n^*}|/n = k_n^*/n$ will be referred to as the \emph{jamming constant}\footnote{This name comes from the vocabulary of parking processes and random sequential adsorption.}. 


\subsection{Greedy algorithm on large sparse random graphs}

The greedy algorithm is the simplest sequential algorithm, and consists in choosing in each step $k$ a new node uniformly at random from ${\mathcal U}_k$. The jamming constant obtained by this algorithm run on a graph $G^{(n)}$ will be denoted by $\sigma_{Gr}(G^{(n)})$.

\subsubsection{Greedy algorithm on sparse Erd{\"o}s-R\'enyi}
\label{sec:greedy_ER}

Let $G^{(n)}$ be a sparse Erd{\"o}s-R\'enyi graph; i.e.\ an $n$-sized graph such that for every pair of vertices an edge exists between them with probability $\lambda/n$ independently of the other edges. The distribution of a graph constructed in this manner will be denoted $ER_n(\lambda/n)$.

The analysis of these graphs is quite simple due to the symmetry and independence of the connections between their vertices. In fact, the process describing the dynamics of the greedy algorithm can be modelled as a one-dimensional Markov process, this being so because after the $k$-th step of the algorithm, if the number of unexplored vertices is $Z^{(n)}_k$, then the unexplored subgraph is still a sparse Erd{\"o}s-R\'enyi graph $G_k \sim ER_{Z^{(n)}_k}(\lambda/Z^{(n)}_k)$. Moreover, in an Erd{\"o}s-R\'enyi graph the degree of vertices are interchangeable (every vertex has, in distribution, the same number of neighbors) and only depends on the parameter $\lambda$ and the size of the graph. 

Using standard fluid limit theorems, it can be shown (see for instance \cite{Bermolen2017}) that:
\begin{gather*}
    \mathbb{E} \left( \sup_{t \in [0,T]} \left| \frac{Z^{(n)}_{\lfloor tn\rfloor }}{n} - z(t) \right| \right) \xrightarrow{n\rightarrow \infty} 0,
\end{gather*}
where $z(t)$ is the solution to the differential equation $\dot z = - \lambda (1 - z)-1$ with the initial condition $z(0)=1$.

Then, it can be proved that in this case $\sigma_{Gr}(G^{(n)})$ converges in $L^1$ to $\tau*$ (defined as the smallest solution to $z(\tau^*)=0$), and thus the jamming constant for large Erd{\"o}s-R\'enyi graphs is arbitrarily close to $\tau^* = \frac{\log (1 + \lambda)}{\lambda}$. See ~\cite{Bermolen2017} where a central limit theorem is also proved.

\subsubsection{Greedy exploration on sparse Configuration Models}
\label{sec:greedy_cm}

The exploration processes of general random graphs cannot be described just by the number of vertices on the unexplored subgraph, as the vertices have degrees that are not interchangeable and that depend in a complicated way on the evolution of the process. This makes the analysis much more involved than in the case of an Erd{\"o}s-R\'enyi graph.

However, a much broader family of random graphs may be considered by resorting to the \emph{Configuration Model}. Given a degree sequence $d^{(n)} = (d_1,...,d_n) \in \mathbb{N}_0^n$, the Configuration Model (see~\cite{wormaldCM, bollobas01, molloy1998size, Durrett, remco2016} and the references therein) is a construction that results in a multi-graph\footnote{A graph with possibly multiple edges between a pair of vertices and self-edges between a vertex and itself.} with the prescribed degrees $d^{(n)}$; i.e.\ one in which the $i$-th vertex will have degree $d_i$ for $i=1,\ldots,n$. To construct such multi-graph, first assign to each vertex a number of \emph{half-edges} equal to its degree. Then, sequentially and randomly pair each half-edge with another unpaired half-edge until there are no more left. The resulting multi-graph will be obtained by assigning an edge between two vertices for each paired couple of half-edges they share. The distribution of a graph constructed according to the Configuration Model will be denoted $CM_n(d^{(n)})$. 

It is known that the order in which the pairing is done does not affect the distribution of the obtained multi-graph~\cite{bollobas01}. Also, although the result of the construction is a multi-graph, standard techniques~\cite{remco2016} can be used to transfer results to simple graphs with the same degrees. 

As shown in~\cite{bermolen2017jamming}, these properties of the Configuration Model allow for the description of the evolution of the greedy algorithm in a Markovian fashion without keeping track of the whole graph structure, but instead only of the degree distribution of unexplored vertices towards other unexplored ones. This is done by constructing simultaneously the random graph and the associated exploration process. As we will see, this makes possible the computation of the resulting jamming constant.

In contrast to the previous discussion, it will now be convenient to describe the algorithm as a continuous time process instead of one taking place in steps. For this, to each unexplored vertex we will associate a random exponential clock of rate $1$. An unexplored vertex will now become active when its clock rings, at which point its associated half-edges are matched and the chosen neighbors which were still unexplored become blocked. Note that this does not affect at all the resulting independent set, since activated nodes are still being chosen uniformly from $\mathcal{U}_k$, as long as all clocks have the same rate. 

Denoting by $\mu^{(n)}_t(i)$ ($i=1,\ldots,n$) the number of unexplored vertices of degree $i$ towards other unexplored vertices (whose initial value $\mu^{(n)}_0(i)$ is simply how many elements in $d^{(n)}$ are equal to $i$), the time evolution of this empirical measure can be proved to be enough to characterize the exploration process in Configuration Model graphs~\cite{bermolen2017jamming}. For this purpose, the associated measure-valued continuous-time Markov process should be scaled according to:

\begin{gather*}
    \bar\mu^{(n)}_t(i)=\frac{1}{n} \mu^{(n)}_t(i),\,t\ge 0,
\end{gather*}
for all $t\geq0$ and $i=1,\ldots,n$.

Then, under the assumption that the initial empirical measure of degrees converges to a measure having mild moment restrictions (i.e.,\ $\bar\mu^{(n)}_0\underset{n}{\rightarrow} \bar\mu_0$), and taking $n$ to infinity, a functional law of large numbers on the evolution of the empirical measure of degrees can be proved~\cite{bermolen2017jamming}. In this work it was also proved that the limit is unique, and given by the solution of a non-linear infinite-dimensional system of differential equations. 

To express this equations in a concise form, denote

\begin{align}
\label{eq:defalpha}
\alpha_t(i)&=\frac{\bar\mu_t(i)}{\sum_{j=0}^{\infty}\bar\mu_t(j)},\\
\label{eq:defbeta}
    \beta_t(i)&=\frac{i\bar\mu_t(i)}{\sum_{j=0}^{\infty} j\bar\mu_t(j)},
\end{align} 
for all $t \geq 0$ and all $i \geq 1$. The probability measures $\alpha_t(\cdot)$ and $\beta_t(\cdot)$ have intuitive interpretations: the first one describes the degree distributions of a randomly (and uniformly) chosen unexplored node at time $t$, while the second one is the \emph{size biased} distribution associated to $\alpha_t(\cdot)$ and represents the degree distribution of a randomly chosen neighbor of a given vertex. This interpretation will be useful when generalizing these results in later sections.

Then, the system of ODE's that describes the evolution of the number of unexplored nodes with degree equal to $i$ towards other unexplored ones can be written as: 

\begin{equation}\label{eq:greedy}
         {\frac{d}{dt}} \bar\mu_t(i) = - \sum_{j=0}^{\infty} \barmu_t(j) \left[ \alpha_t(i) + F_t(i) \right],
\end{equation}

where $F_t(i):=\sum_{k=0}^{\infty} k\alpha_t(k)\left(\beta_t(i+1) + (\beta_t(i)-\beta_t(i+1))\sum_{l=0}^{\infty} l\beta_t(l)\right)$. This equations also have a simple and intuitive interpretation. The first sum on the right hand side is the total transition rate of the Markov process. The first term within the square brackets represents the possibility that a node of degree $i$ is activated at time $t$, in which case the number of unexplored nodes with degree equal to $i$ decreases by one. While the term $F_t(i)$ corresponds to the blocking of neighbors of the selected node and the removal of their edges from the unexplored subgraph. Note that the activated node has degree distribution $\alpha_t(\cdot)$, while the blocked neighbors and their neighbors all have degree distribution $\beta_t(\cdot)$. This representation will be useful to analyze the degree-aware algorithms in Sec.\ \ref{sec:awaredynamic}.

The main consequence of this limit is the characterization of the jamming constant of the greedy algorithm run on Configuration Models~\cite{bermolen2017jamming}:

\begin{thm}\label{theo.spa} 
For each $n\geq 1$, let $G^{(n)}\sim CM_n(d^{(n)})$. Under certain moment assumptions (see Theorem 3.1 in \cite{bermolen2017jamming}), the following result holds:
 \begin{equation*}
     \mathbb{E}(|\sigma_{Gr}(G^{(n)}) - c_{\bar\mu_0} |) \xrightarrow[n \rightarrow \infty]{} 0, 
 \end{equation*}
 where
 \begin{equation}\label{eq:int_greedy}
      c_{\bar\mu_0} = \int_0^\infty \sum_{i=0}^{\infty} \bar \mu_t(i) dt. 
 \end{equation}
\end{thm}

A more explicit result for the jamming constant was later proved in \cite{janson2017greedy} for a modified dynamics that allows for a simplification of the limiting differential equation system. In this work, the authors studied a different hydrodynamic limit which results in simpler equations that can be directly integrated:

\begin{thm}
Under a second moment assumption on the initial distribution $\bar\mu_0$, and
calling $\tau_\infty$ the unique value in $(0, \infty]$ such that
 \begin{equation*}
     \int_0^{\tau_\infty} \frac{e^{-2 h}}{\sum_{i=0}^{\infty} i \bar\mu_0(i) e^{-i h}} dh = 1, 
 \end{equation*}
 then
 \begin{equation}\label{eq:int_greedy_janson}
     \sigma_{Gr}(G^{(n)}) \xrightarrow[n\rightarrow \infty]{\mathbb{P}} \int_0^{\tau_\infty} \frac{e^{-2h}\sum\limits_{i=0}^{\infty} \bar\mu_0(i) e^{-i h}}{\sum\limits_{i=0}^{\infty} i \bar\mu_0(i) e^{-i h}} dh. 
 \end{equation}
\end{thm}
 

\subsection{Degree-greedy algorithm on regular graphs}
The degree-greedy algorithm is a variant of the greedy that takes into account the degree of the vertices in the unexplored subgraph. The algorithm is exactly as the one described in the previous section, except that at the $k$-th step a vertex $v$ is selected uniformly from the vertices of \emph{minimum degree} within the subgraph of unexplored vertices $G_k$. The rest of the exploration is as before: the state of the chosen vertex is changed to active, and its neighbors to blocked. The algorithm ends when there are no more vertices in the unexplored subgraph. We will denote by $\sigma_{DGr}(G^{(n)})$ the jamming constant obtained by the degree-greedy algorithm ran on $G^{(n)}$. 

Although it has been studied in the computer science community (for example, in \cite{halldorsson1997greed}), there are very few mathematical results that characterize or bound the independent set found by this algorithm. The main result on this respect was presented by Wormald in \cite{wormald1995differential}. In this work, a fluid limit for the process generated by the degree-greedy algorithm when run on a $d$-regular graph (a graph constructed by the Configuration Model where every vertex has the same degree $d\geq1$) was proved. Numerical estimations based on this fluid limit are provided for different values of $d$ in \cite{wormald1995differential}. Also, in \cite{karp1981maximum} the authors give a characterization of a matching constructing degree-greedy variant, that under certain circumstances can be used to determine the behaviour of the degree-greedy algorithm ran on Erd{\"o}s-R\'enyi graphs. 


\subsection{Dynamical exploration}\label{sec:glauber}

As opposed to the sequential character of previous algorithms, more complex dynamics where nodes join and leave the independent set can be considered. The simplest dynamics of this type would be the Glauber dynamics where each node independently tries to join the set of active nodes (and succeeds in the absence of interfering activated nodes) but also leaves this set (\emph{deactivates}) after a random time. As we will see, Glauber dynamics allow to define a limiting stochastic process that approaches a maximum independent set. Moreover, as we further discuss in the simulations section, this dynamics is very similar to that of 802.11-based nodes. 

Given a graph $G^{(n)}$, the Glauber dynamics is defined as a discrete time Markov process on the state space of subsets of vertices $S = \mathcal{P}(V)$. At each step $k \in \mathbb{N}_0$, the set $A_k \subseteq V$ will be interpreted as the set of vertices \emph{transmitting information} at that time (that is, the active vertices) and will be referred as a \emph{configuration}. Given a configuration $A_k$, the configuration at the following step $A_{k+1}$ will be constructed  according to the following rules:
\begin{itemize}
 \item Choose a vertex $v$ uniformly from $V$.
 \item With fixed probability $\frac{\beta}{1+\beta} > 0$, if no neighbor of $v$ belongs to $A_k$, set $A_{k+1}= A_k \cup \{v\}$.
 \item With probability $\frac{1}{1+\beta} > 0$, set $A_{k+1} = A_k \backslash \{v\}$.
\end{itemize}
It is easy to see that if the initial configuration $A_0$ is an independent set, $A_k$ will be an independent set for every time $k \in \mathbb{N}_0$. The invariant measure is known to be given by (for every configuration $A \subseteq V$)~\cite{vigoda2001note}:
\begin{equation}
	\mu( A ) = \frac{\beta^{|A|}}{Z},
\end{equation}
where $Z$ is a normalization constant. This means that in the limit $\beta \rightarrow \infty$, the invariant measure concentrates in the maximum independent sets. In some special situations, the mixing time has been characterized. For example in \cite{vigoda2001note}, it was proved that when the degree distribution is bounded by $\Delta \geq 0$ and $\beta < \frac{2}{\Delta - 2}$ the mixing time is $\mathcal{O}(n \log(n))$. Or in the case of a bipartite regular graph, it was shown \cite{galvin2006slow} that for $\beta$ large enough the mixing time is exponential in $n$.

In the limit of $\beta \rightarrow \infty$, Glauber dynamics will behave in the following manner: 
\begin{itemize}
    \item In a first phase (taking $A_0$ as the empty set), vertices will activate until a maximal independent set is reached. The distribution of the size of the resulting independent set will match that of the greedy algorithm, since nodes are chosen uniformly. 
    \item Because activation attempts occur much faster than desactivations, after this initial phase there will be plenty of failed activation attempts (which can be omitted in the analysis) followed by a single deactivation of some node $v$. 
    \item After this single deactivation happens, the vertices in $\mathcal{N}(v) \cup \{ v \}$ that do not have active neighbors will try to activate uniformly until no more activations are possible and a new maximal independent set is reached.
    \item The last two steps are repeated alternatively.
\end{itemize}

Because the invariant measure in this limit concentrates on independent sets of maximum size, this dynamic will asymptotically approach a maximum independent set, but this can take a prohibitively long time. In Sec.\ \ref{sec:maxywifi}, we will discuss the application of this kind of limit dynamics to communication networks modeled by the Configuration Model.


\subsection{Characterization of maximum independent sets}

Finally, we review some of the results that characterize or bound the asymptotic independence number of the Erd{\"o}s-R\'enyi graphs and the Configuration Model. In the case of the Configuration Model, we will focus on the $d$-regular graphs as most literature centers around them. As we will see, many of these bounds are estimations that result from the analysis of independent set finding algorithms.

\subsubsection{Erd{\"o}s-R\'enyi random graphs}

In the case of sparse Erd{\"o}s-R\'enyi graphs, in the general case, exact values or even upper bounds on the independence number are unknown. A lower bound is provided by the greedy algorithm. Finer results in that direction are proved in \cite{bollobas1976cliques,frieze1990independence}:
\begin{thm}
    For each $n\geq1$ and $\lambda> 3$,  let $G^{(n)} \sim \mbox{ER}_n(\lambda/n)$, then the independence number is bounded by:
  $$ \alpha (G^{(n)}) \leq \frac{2\log(\lambda)}{\log(1-\lambda/n)}$$
\end{thm}
Since the size of the independent set discovered by the greedy algorithm $n\sigma_{Gr}(G^{(n)})$ is known asymptotially (see Section ~\ref{sec:greedy_ER}) it can be proved that:

$${\mathbb P}\Big(n\sigma_{Gr}(G^{(n)}) \geq (1+\epsilon)\frac{\alpha (G^{(n)})}{2}\Big) \underset{n\rightarrow \infty}{\rightarrow} 1,$$
that is, the greedy algorithm discovers independent sets whose size are asymptotically larger than half of the independence number.

Furthermore, in \cite{karp1981maximum} the authors study the problem of constructing maximum matchings; their results imply that the degree-greedy algorithm is asymptotically optimal for Erd{\"o}s-R\'enyi graphs when the mean degree is $\lambda < e$ and determine the value of the independence number in this case. To reach this conclusion, some optimality lemmas similar to the ones we prove in later sections (but in their case for matchings) are stated. However, our result applies to a broader family of random graphs.
\\

\subsubsection{Regular random graphs}

The proof of the existence of a limiting independence number for random $d$-regular graphs was given in \cite{bayati2010combinatorial}. Moreover, in \cite{lauer2007large} Lauer and Wormald proved that the independence number of a $d$-regular graph is w.h.p.\ bounded from below by 
\begin{equation*}
    \beta(d) =\frac{1}{2}\left[1-(d-1)^{-2/(d-1)}\right].
\end{equation*}

In the same work, Wormald (and Gamarnik and Goldberg independently in~\cite{gamarnik2010randomized}) showed that, when the girth of the graph (the length of its shortest cycle) goes to infinity, the proportion of vertices in the independent set found by a greedy algorithm in a $d$-regular graph is given by $\beta(d)$ (for $d\geq3$), asymptotically (in probability) when $n\rightarrow\infty$. These results can also be derived by means of the theorems presented in Sec.\ \ref{sec:greedy_cm}.

While upper bounds were also proved in \cite{bollobas1981independence,mckay1987lnl} and an alternative lower bound in \cite{frieze1990independence}. In a recent work~\cite{ding2016maximum}, the exact law of large number for the independence ratio of $d$-regular graphs of sufficiently large $d$ was established. 

In Sec. \ref{sec:analysisapproximation}, we will compare some of these bounds with the approximate performance of the degree-greedy algorithm.

\subsection{Optimality of the degree greedy algorithm}\label{sec:optimality}

In this section, we detail sufficient conditions for the degree-greedy algorithm to find an independent set which is asymptotically of maximum size, as the number of nodes grows.
The proof (which are mathematically involved and do not fall in the scope of this paper) can be found in a companion paper dedicated to this result.

A sufficient condition for the degree-greedy algorithm to find w.h.p.\ an independent set that asymptotically contains the same proportion of vertices as a maximum one can be proved: 

\begin{proposition}\label{prop:condnearopt}[\cite{jonckheere2018asymptotic}]
    Given a sequence of graphs $G^{(n)}$ distributed according to $\mbox{CM}_n( d^{(n)})$, assume that:
\begin{itemize}    
  \item In a first phase, the degree-greedy algorithm defines w.h.p.\ a selection sequence that selects only vertices of degree 1 or 0. Let $T$ be the limiting length of this phase and $\mu_T$ the limiting degree distribution of the graph at the end of this phase. 
  \item
  Assume further that
  $$ \frac{\sum_{i} i(i-1) \mu_T(i)}{\sum_{i} i \mu_T(i)} < 1,$$ 
  i.e., at the end of this first phase, the remaining graph is subcritical,
\end{itemize}   
    then
     $$  \frac{\sigma_{DGr}(G^{(n)})}{\alpha(G^{(n)})}  \xrightarrow[n \rightarrow \infty]{} 1, \text{ in probability}.$$
\end{proposition}    

Where $\sigma_{DGr}(G^{(n)})$ represents the (random) proportion of vertices in the independent set found by the Degree-Greedy. The intuition behind this result is that a subcritical graph \emph{does not differ much} from a collection of trees. We can then couple the algorithm running in the subcritical graph with one running in the collection of spanning trees of its components. The difference of the independent sets obtained by both coupled algorithms will be of negligible size as a subcritical graph has few components that are not trees and the degree-greedy algorithm running in the trees will select only vertices of degree 1 or 0.

The importance of this proposition lies in the fact that, when the graph gets smaller, the algorithm may very well choose from time to time vertices with degrees larger than 1. What proposition 3.4 shows is that this does not spoil asymptotical optimality if the selection of this vertices happen after the graph has broken into a collection of small components (i.e., when the graph has already become subcritical).

\subsection{Practical consequences}

\subsubsection{Verifying optimality}

As a direct consequence of the results of previous section, the degree-greedy exploration on any subcritical graph is optimal. In other words, if the limiting degree distribution of the configuration model is such that $\frac{\sum_{i} i(i-1) \mu_0(i)}{\sum_{i} i \mu_0(i)} < 1$, then the degree-greedy strategy is asymptotically optimal.

However, there is much wider class of graphs for which this property holds as this result also deals with super-critical graphs such that the degree-greedy algorithm might select only degree 1 or 0 for a sufficient long first phase. We can actually  numerically characterize for a given (parametrized) distribution the set of parameters such that the degree-greedy algorithm is optimal. It indeed suffices to numerically estimate the hydrodynamic limit and show that during a first phase, only degree 0 or 1 are explored while the resulting graph after this first phase is sub-critical.

We mention here two examples. Consider first graphs with asymptotic degree distribution given by $\mbox{Pois}(\lambda)$. Then the subcriticality condition is given by the condition $\lambda < 1$. Estimations for the hydrodynamic limit of the degree-greedy dynamics were made using groups of 10 Configuration Model graphs of 50.000 vertices each with Poisson distributions of varying mean degrees that ranged from 1 to 2.5. In all the simulations it was consistently found that conditions of proposition 3.4 are met. This means that the degree-greedy algorithm remains asymptotically optimal for graphs of mean degree smaller than 2.5, way after the critical connectivity threshold. This is consistent with the results proved in \cite{karp1981maximum} and \cite{jonckheere2018asymptotic} where it is proved that the optimality threshold is $e$ (exponential number).
This is also in accordance with our simulations of the next section (see Figure \ref{fig:comp_cotas_er}). 
 
Another example is a power-law distribution. We can show numerically that when the exponent of the power law satifies $a > 3$, then the optimality result holds. It is of course interesting to provide theoretical bounds or charaterizations for these thresholds depending on the original degree-distribution but this is out of the scope of the present work. For further details on these kind of proofs, see \cite{jonckheere2018asymptotic}.

\subsubsection{Heuristics for efficient discovering algorithms}
Since the degree-greedy algorithm works (quasi)-optimally for graphs with low enough connectivity, and that Glauber dynamics are theoretically optimal in the long run but get stuck in practice in ``local minima'' configurations, one can unite the best of both worlds and to start Glauber dynamics with a configuration found previously by a degree-greedy strategy. We further explore this approach in Sec.\ \ref{sec:maxywifi}.


\section{Generalization of sequential algorithms}\label{sec:generalization_seq_algorithms}

In this section we adapt some of the results reviewed in Sec.\ \ref{sec:greedy_cm} for the greedy algorithm on the Configuration Model in order to obtain more general results that allow for a degree-aware control. Through these results, we are able to propose decentralized schemes
with a performance comparing very closely to the one of degree-greedy.

In Sec.\ \ref{sec:awaredynamic}, we analyze the dynamical case. As we shall see, the degree-greedy algorithm can be thought of as a limit process for this family of dynamics, allowing this processes to be used to approximate it. On the other hand, in Sec.\ \ref{sec:awarestatic} we will analyze the static case which results in a simpler asymptotic behaviour that can be described by solving just a single ODE, allowing for easier simulations and implementation.


\subsection{Dynamic degree-aware exploration}\label{sec:awaredynamic}

In section \ref{sec:greedy_cm}, the greedy dynamics was presented along with a set of differential equations that describe its limiting behaviour. Here we consider a generalization of this in which each vertex with degree $i \in \mathbb{N}_0$ towards other unexplored vertices has an exponential clock of parameter $0<\lambda(i)<\infty$. The blocking of the neighbors of the selected vertices, and the matching of the half-edges of these neighbors towards the set of unexplored vertices are done exactly as before (that is, uniformly at random). We will call this process the \emph{dynamic degree-aware} exploration.

Observe that changing the rates of these exponential clocks is equivalent to changing the probability that in a transition at time $t\geq0$ a vertex of degree $i \in \mathbb{N}_0$ is activated: from $\alpha_t(i)$ (which corresponds to the degree algorithm, defined in section \ref{sec:greedy_cm}) to 
\begin{equation}\label{eq:nuevaproba}
    \gamma_t(i) = \frac{\lambda(i) \bar \mu^{(n)}_t(i) }{\sum\limits_{j=0}^{\infty} \lambda(j) \bar \mu^{(n)}_t(j)},
\end{equation}
where $\{ \bar \mu^{(n)}_t(i)\}_{i\geq0}$ is (as in Sec.\ \ref{sec:greedy_cm}) the scaled number of unexplored vertices with degree $i$ towards other unexplored vertices.

Then, for the process associated with $\{ \bar \mu^{(n)}_t(i)\}$ the following limit holds:

\begin{proposition}
Under the same assumptions of Theorem \ref{theo.spa} and letting $(\lambda(i))_{i\geq0}$ be bounded and positive, then the sequence of processes $\{\bar\mu^{(n)}_t\}$ converges in probability and uniformly on compact time intervals towards the only measure-valued deterministic function $\bar \mu_t$ solution of the following infinite dimensional differential system: 
\begin{equation}
    \frac{d}{dt} u_t(i) = -\sum_{j=0}^{\infty} \lambda(j)u_t(j)\left[\gamma_t(i) + G_t(i)\right], \label{eqdiffdificil}
\end{equation}
where $G_t(i):=\sum_{k=0}^{\infty} k\gamma_t(k)\left(\beta_t(i+1) + (\beta_t(i)-\beta_t(i+1))\sum_{l=0}^{\infty} l\beta_t(l)\right)$ with $\gamma_t(\cdot)$ given by \eqref{eq:nuevaproba} and $\beta_t(\cdot)$ defined as before in \eqref{eq:defbeta} (both functions associated to $(u_t(i))_{i\geq0}$).
\end{proposition}
\noindent
\emph{Proof.} The main difference in the processes in question is that, for the greedy algorithm, an unexplored vertex is activated when its exponential clock of rate $1$ rings; while in these, the rates of the clocks are not necessarily the same and can depend on their degrees: initial degrees (in the \emph{static} models) and actual degrees towards other unexplored vertices (in the \emph{dynamic} ones). 

Note that these equations are just adaptations of \eqref{eq:greedy} where we substitute $\alpha_t(i)$ by $\gamma_t(i)$ and the total transition rate is now $\sum_{j=0}^{\infty} \lambda(j)u_t(j)$. The rest of the equations remains unchanged since, as explained before, the neighbors of the active vertices are in both cases selected uniformly at random. The proof is then, mutatis mutandi, the same as the proof of Theorem \ref{theo.spa} of \cite{bermolen2017jamming} which is quite involved and we do not reproduce here. The interested reader is referred to \cite{bermolen2017jamming}. $\square$

Once equations \eqref{eqdiffdificil} are solved, the corresponding jamming constant can be computed according to:
\begin{gather}
    c_{\bar\mu_0} = \int_0^\infty \sum_{i=0}^{\infty} \lambda(i)\bar\mu_t(i)dt,\label{eqsdificil}
\end{gather}
where $\bar \mu_0(\cdot)$ is the scaled limit of the initial degree distribution. As before, if we denote by $\sigma_{DDA}(G^{(n)})$ the proportion of vertices in the independent set found, it will converge in $L^1$ to $c_{\bar\mu_0}$ when $n$ goes to infinity.

The degree-greedy algorithm could be considered in this context by setting the probabilities $\gamma_t(i)$ to be given by
\begin{gather}\label{greedyposta}
    \gamma^*_t(i) = 
    \begin{cases}
        1, \quad \text{if }\bar\mu_t(i)>0 \text{ and } \bar\mu_t(j)=0\; \forall j<i;\\
        0, \quad else.
    \end{cases}
\end{gather}

However, among other difficulties, these probabilities are not of the form \eqref{eq:nuevaproba} and therefore the limit described by the previous theorem (equations \eqref{eqdiffdificil}) does not apply. Nonetheless, we can consider families of clock rates $\lambda(\cdot)$ that will select at each transition nodes with minimum degree with high probability. For example, $\lambda(i) = (i+1)^{-L}$ (with $L > 0$). In this case, the associated probabilities will be given by
\begin{equation*}
    \gamma_t(i) = \frac{(i+1)^{-L} \bar \mu_t(i)}{\sum\limits_{i=0}^{\infty} (i+1)^{-L} \bar \mu_t(i)}.
\end{equation*}

Note that these $\gamma_t(i)$ tend to $\gamma_t^*(i)$ as $L$ goes to infinity, and so this family of processes gives a way of approximating the size of the independent set discovered by the degree-greedy algorithm. In Sec. ~\ref{sec:simulaciones} we will compare this approximation with the results obtained by using the probabilities $\gamma^*_t(\cdot)$. As we will see, the proposed approximation presents an excellent performance for different initial degree distributions.

As a consequence of Proposition \ref{prop:condnearopt}, this approximation of the degree-greedy algorithm can be used to estimate the independence number of a wide family of random graphs, which is in general a hard problem. Moreover, note that the algorithm may be implemented in a distributed manner among vertices, as they only need to know their degree towards unexplored nodes. 

\begin{remark}
Nor the techniques neither the results of \cite{janson2017greedy}, which will be of use in the next section, can simplify the limit presented in this section. Indeed, in the processes under study the information of the degree distribution of unexplored vertices towards other unexplored is needed at each time to derive limit equations, which renders the strategy in \cite{janson2017greedy} inapplicable.
\end{remark}


\subsection{Static degree-aware exploration}\label{sec:awarestatic}

Here we will discuss another variation of the greedy dynamics which leads to a simpler limit. Similarly to the situation discussed in the previous section, each unexplored vertex will activate when an exponential clock rings but instead of every vertex having a clock that dynamically depends on its degree in the unexplored subgraph, now its rate will only depend on the initial degree of the vertex in question and will therefore be constant in time. 

Again, the blocking of vertices and matching of half-edges is done as in the regular greedy algorithm. We will call this process the \emph{static degree-aware} exploration, and we will denote the proportion of the vertices in the independent set found by it in a graph $G^{(n)}$ by $\sigma_{SDA}(G^{(n)})$. In this context, the following result holds:

\begin{proposition}\label{prop:staticDA}
For $n\geq1$, let $G^{(n)}\sim \mbox{CM}_n(d^{(n)})$ with asymptotic degree distribution $(\bar\mu_0(k))_{k\in\mathbb{N}_0}$ (of mean $m := \sum_{k=0}^\infty k \bar\mu_0(k)$) and assume the convergence of the second moment of the degree distribution towards $\sum_{k=0}^\infty k^2 \bar\mu_0(k)$. Then, if $(\lambda(k))_{k\in\mathbb{N}_0}$ is such that $\sum_{k=0}^\infty k \lambda(k) \bar\mu_0(k)$ is uniformly summable, we will have that

\begin{equation*}
    \sigma_{SDA}(G^{(n)}) \xrightarrow[\mathbb{P}]{n\rightarrow\infty} \int_{0}^\infty \sum_{k=0}^{\infty} \lambda(k) \bar \mu_0(k) e^{-\lambda(k) t} e^{- k \tau(t)} dt,
\end{equation*}

where $\tau(t)$ is the solution to the ODE 

\begin{equation*}
    \frac{d \tau}{dt} = \frac{1}{m}\sum_{k=1}^{\infty} k \lambda(k) \bar\mu_0(k) e^{-\lambda(k) t} e^{-(k-2) \tau}
\end{equation*}

with initial condition given by $\tau(t)=0$.

\end{proposition}
\noindent
\emph{Proof.} The proof of the convergence of the dynamics is a modification of the one in \cite{janson2017greedy}.

First, the process is described as a Markov process that involves the coordinates of the degree distribution of unexplored vertices ($\bar\mu_t(k)$, for $k\geq0$), the total number of unmatched half-edges ($U_t$) and the number of active vertices ($A_t$).

Second, the drifts $\delta(\cdot)$ for all the coordinates of the process are explicitly obtained. In our case, they will be 

\begin{equation*}
    \delta(A_t)/n= \sum_{k=0}^\infty \lambda(k) \bar\mu_t(k),
\end{equation*}

\begin{equation*}
    \delta(U_t)/n = - \sum_{k=0}^\infty k \lambda(k) \bar\mu_t(k) \left( 2 - \frac{k-1}{U_t-1}\right),
\end{equation*}

\begin{equation*}
    \mbox{and }\delta(\bar\mu_t(k)) = - \lambda(k) \bar\mu_t(k) - \sum_{j=0}^\infty n p_{jk}\lambda(j) \bar\mu_t(j) (\bar\mu_t(k)-\delta_{jk}/n),
\end{equation*}
respectively. Here, $p_{jk}$ represents the probability of connection between a pair of nodes of degree $j$ and $k$.

By Dynkin's formula, the forms of the drifts and bounding the corresponding quadratic variation of the martingales resulting, the convergence of each of these hydrodynamic limits are established for some subsubsequence of every subsequence.

After this, the asymptotic form of the drifts need to be obtained. The bounds and asymptotic forms are analogous to the ones in \cite{janson2017greedy}. 
    
However, to prove $\delta(U_t)/n \rightarrow -2 \sum_{k=0}^\infty k \lambda(k)\bar\mu_t(k)$ we need a different argument to show that the number of unmatched half-edges diverges at any given time as $n\rightarrow\infty$, in particular if we want to assume that $\lambda(k)$ is unbounded (The bounded case can be handled as in \cite{janson2017greedy}).
    
To prove this, note that because $\sum_{k=1}^\infty k \lambda(k) \bar \mu_t(k)$ is uniformly summable, there exists a sequence $(\epsilon_N)_{N\geq1}$ s.t. (for every $n\geq1$) $\sum_N^\infty k \lambda(k) \bar \mu_t(k) \leq \epsilon_N$ and $\epsilon_N \xrightarrow{N\rightarrow\infty}0$. Then using Dynkin's formula,
    
    \begin{equation}\label{eq:cotaU}
        U_t/n \geq U_0/n - \int_0^t \sum_{k=1}^\infty k \lambda(k) \bar \mu_s(k) ds \geq U_0/n - \bar\lambda_N \int_0^t \sum_{k=1}^N k \bar \mu_s(k) ds - \epsilon_N t
    \end{equation}

where $\bar\lambda_{N} := \max_{k\in\{1,...,N\}} \lambda(k)$. Calling $u_t$ the limit of $U_t/n$ in the corresponding subsubsequence, because of \eqref{eq:cotaU} and the fact that $U_t/n \geq \sum_{k=1}^\infty k \bar\mu_t(k)$, we have that for every $t\geq0$ and $\bar N\geq1$ large enough $u_t \geq \bar u_t$, where $\bar u_t$ is the solution of the differential equation $\bar u_t' = - \bar\lambda_{\bar N} \bar u_t - \epsilon_{\bar N}$ with initial condition $\bar u_0 = m$. This can be integrated to give $\bar u_t = \left( m + \epsilon_{\bar N} / \bar \lambda_{\bar N} \right) e^{-\bar\lambda_{\bar N} t} - \epsilon_{\bar N} / \bar \lambda_{\bar N}$. This further implies that $u_t \geq \bar u_t > 0$ for every $t \in [0,\bar t_{\bar N})$ with
    
\begin{equation*}
        \bar t_{\bar N} = \log \left[ \left( 1 + \frac{m \bar \lambda_{\bar N}}{\epsilon_{\bar N}} \right)^{1/\bar\lambda_{\bar N}} \right] = \frac{m}{\epsilon_{\bar N}} +o(1)
\end{equation*}
    
Note that we used that $\bar\lambda_{\bar N}\xrightarrow{\bar N\rightarrow\infty}\infty$ as $\lambda(k)$ is unbounded. Because $\epsilon_{\bar N}\xrightarrow{\bar N \rightarrow\infty} 0$, we will have that $\bar t_{\bar N} \xrightarrow{\bar N\rightarrow\infty}\infty$. This finally shows that for every fixed $t\geq0$, $u_t > 0$ which in turn means that $U_t \xrightarrow{n\rightarrow\infty}\infty$.
    
Finally, by topological arguments, this convergence that is valid for these subsubsequences is extended to the whole sequence of processes.
 
Then, the hydrodynamic limit of the degree distribution of unexplored vertices in a static degree-aware dynamics will be given by (for $k \geq 1$)
\begin{equation}\label{eq:fluidmodulated}
    \bar \mu_t(k) = \bar\mu_0(k) e^{-\lambda(k) t} e^{- k \tau(t)}, 
\end{equation}
where $\tau(t)$ is as in the statement of the proposition. Because of the convergence of the degree distribution and the number of active vertices, we will have that

\begin{align}
     \sigma_{SDA}(G^{(n)}) & = \lim_{t\rightarrow\infty} \lim_{n\rightarrow\infty} \frac{A_t}{n} = \int_{0}^\infty \sum_{k=0}^{\infty} \lambda(k) \bar \mu_t(k) dt  \\
     & =  \int_{0}^\infty \sum_{k=0}^{\infty} \lambda(k) \bar \mu_0(k) e^{-\lambda(k) t} e^{- k \tau(t)} dt,
     \label{eq:modulatedIS}
\end{align}

plus a term that vanishes in probability. $\square$

Note that this proposition reduces the problem of analyzing the new dynamics to integrating a single ODE, which is much simpler than solving the system presented in the previous section, although the resulting independent set will be naturally smaller than in the dynamic degree-aware exploration case. 

Moreover, calling $c_{\bar\mu_0}$ the limiting jamming constant obtained, the proportion of vertices of degree $i\geq1$ in the independent set will be asymptotically given by
\begin{equation}\label{eq:modulateddegree}
    q_i = \frac{1}{c_{\bar\mu_0}} \int_{0}^\infty \lambda(i) \bar\mu_0(i) e^{-\lambda(i) t} e^{- i \tau(t)} dt.
\end{equation}

Proposition \ref{prop:staticDA} and equation \eqref{eq:modulateddegree} will be used in Sec.\ \ref{sec:fairness} to analyze the  degree distribution of the nodes belonging to the independent set and to derive strategies to improve the equality of access probability (i.e., fairness).


\section{Analysis of the approximation of the degree-greedy algorithm}\label{sec:analysisapproximation}
\label{sec:simulaciones}


In this section we will numerically study the approximations of the degree-greedy algorithm presented in Secs.\ \ref{sec:awaredynamic} and \ref{sec:awarestatic}. Unless stated otherwise, simulations will correspond to 10 random graphs $G^{(n)}$, with $n=1000$, generated according to the Configuration Model corresponding to the considered scenario. For each of these graphs we will compute its empirical degree measure $\mu(\cdot)$, numerically solve \eqref{eqdiffdificil} or \eqref{eq:fluidmodulated} (corresponding to the dynamic and static degree aware algorithms respectively) with clock parameters $\lambda(i) = (i+1)^{-L}$ (with $L>0$ sufficiently large) using as initial condition $\bar\mu_0(\cdot)=\mu(\cdot)/n$, and calculate the size of the maximal independent set through \eqref{eqsdificil} or \eqref{eq:modulatedIS}. 

We will also execute a single degree-greedy exploration process on each of these 10 graphs and report the size of the resulting maximal independent sets in the form of a boxplot. Moreover, we also execute a long enough Glauber dynamics so that it should get close to the maximum size. However, as noted before, this dynamics provides asymptotically the independence number, and it might not be reached through a finite-time simulation, as the mixing time can be exponential in $n$ depending on the graph characteristics. This aspect will be further discussed in Sec.\ \ref{sec:maxywifi}. 

Finally, and for contrasting, we will also include the size of the maximal independent set corresponding to the greedy algorithm (which may be obtained by any of the methods presented in \cite{bermolen2017jamming,janson2017greedy}). 

\subsection{Random regular graphs}

First, we will consider $d$-regular random graphs. As we discussed in Sec.\ \ref{sec:state_of_the_art}, this is the most studied case in the literature. We will also include in our comparison the lower bound on the independence number proved in \cite{wormald1995differential} (see the numerical results presented in Table 1 in that article), and the recent exact value presented in \cite{ding2016maximum}. 

\begin{figure}
    \centering
    \includegraphics[width=0.6\textwidth]{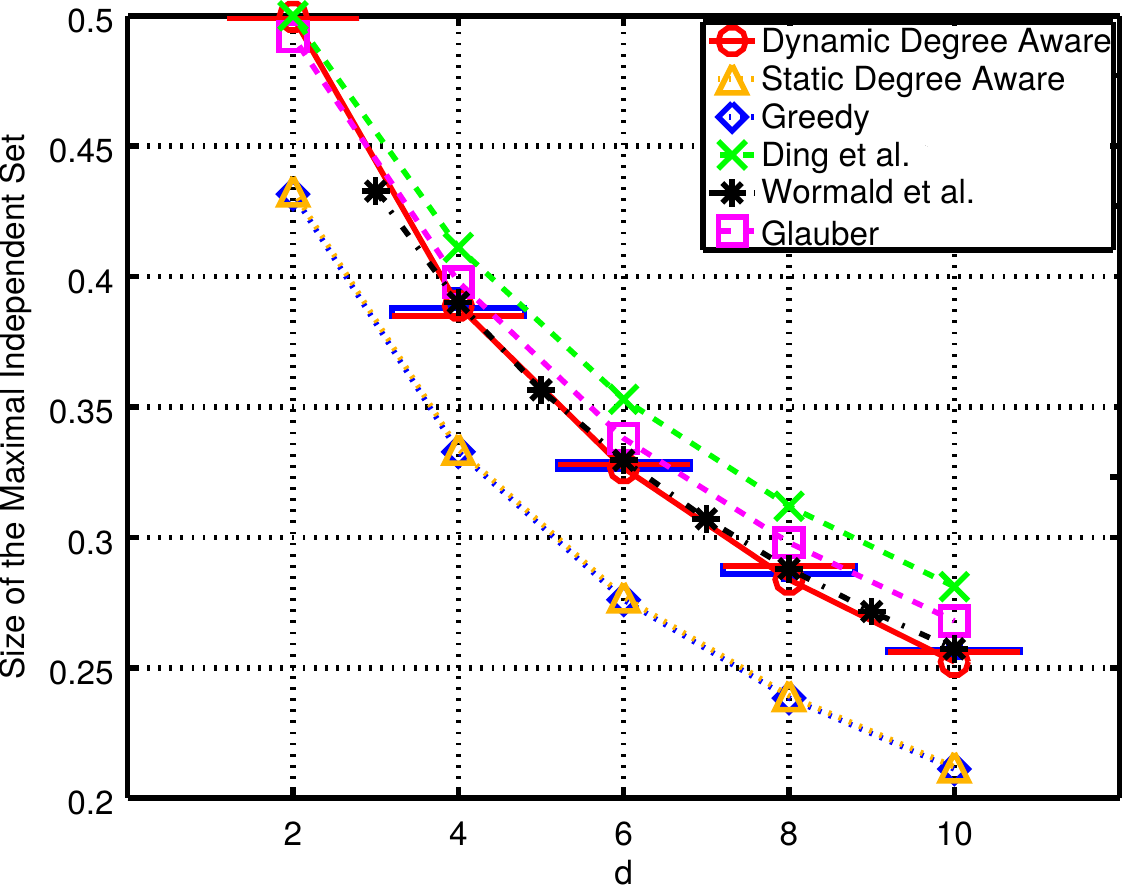}
    \caption{Size of the maximal independent set for different values of $d$ for a $d$-regular random graph. Simulation results for the degree-greedy are shown as a boxplot (seen as horizontal lines), red circles correspond to \eqref{eqsdificil}, yellow triangles to \eqref{eq:modulatedIS}, blue diamonds to \eqref{eq:int_greedy} or \eqref{eq:int_greedy_janson}, green crosses to the values in \cite{ding2016maximum}, black stars to the lower bounds in \cite{wormald1995differential} and magenta squares to the Glauber dynamics.}\label{fig:comp_cotas_d_regular}
\end{figure}

The results are shown in Fig.\ \ref{fig:comp_cotas_d_regular}. Note that, as expected in this case, our analysis of the dynamic degree-aware exploration provides a very accurate approximation of the size of independent set obtained by the degree-greedy algorithm. Moreover, the lower bounds provided by Wormald in \cite{wormald1995differential} give essentially the same values as our approximation. This is due to the fact that, as suggested by Wormald~\cite{wormaldCM}, this lower bound is likely the exact limiting value. 

Furthermore, by observing the results corresponding to \cite{ding2016maximum}, we may conclude that there is possibly still room for increasing the size of the independent set of the degree-greedy algorithm. However, it must be noted that this result is exact only for $d$ larger than a $d_0$ not specified by the authors. Our Glauber dynamics simulations further support that, at least for these values of $d$, those results seem to be upperbounds. In any case, this improvement would be relatively modest, and certainly smaller than the gap between the degree-greedy and the greedy algorithms. Finally, note that in this case (where all nodes have the same initial degree), the static degree aware and the degree algorithms are equivalent and thus obtain the same result.

\subsection{Erd{\"o}s-R\'enyi graphs}\label{sec:erdossimu1}

\begin{figure}
    \centering
    \includegraphics[width=0.6\textwidth]{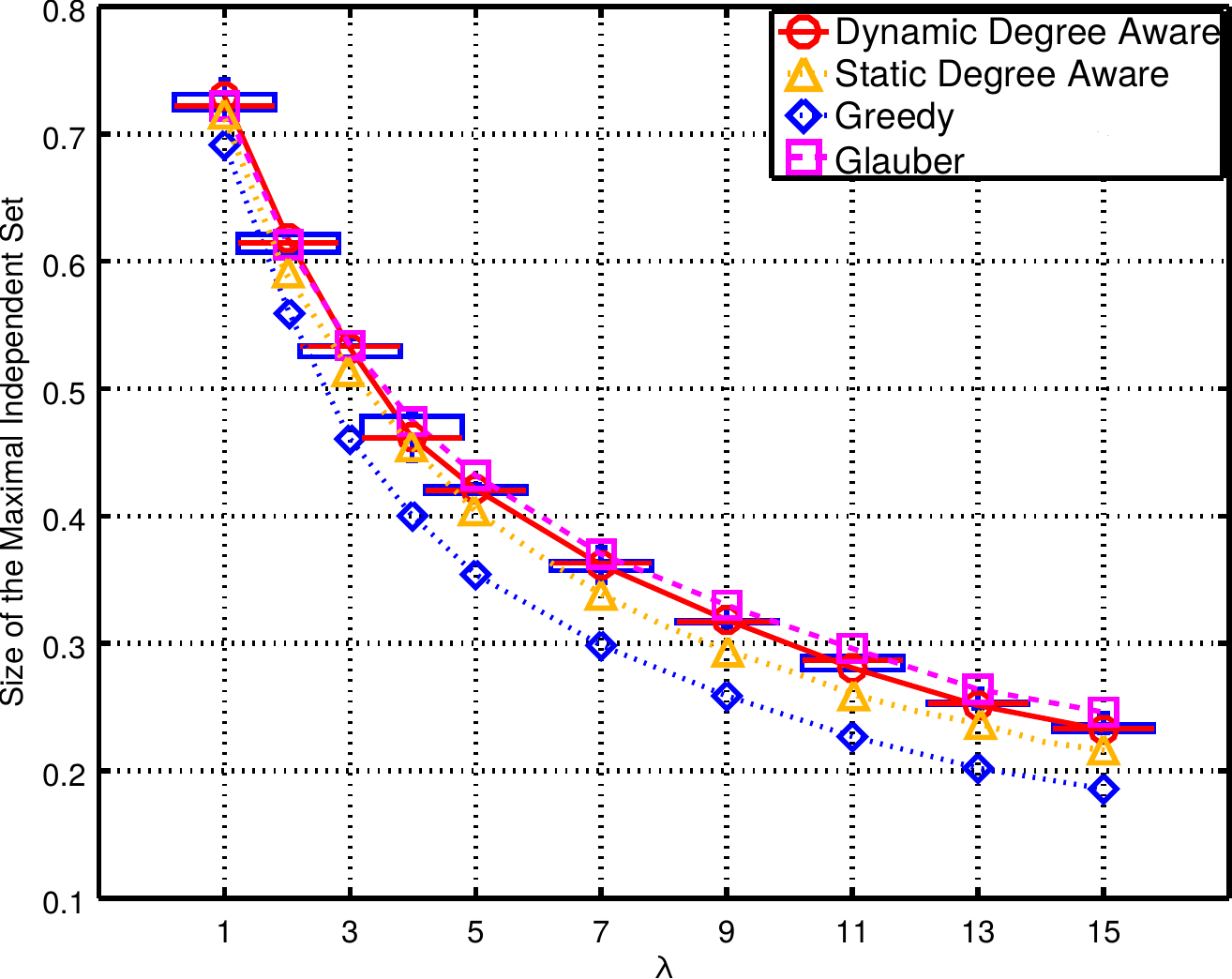}
    \caption{Size of the maximal independent set for different values of $\lambda$ for an Erd{\"o}s-R\'enyi graph with a connecting probability equal to $\lambda/n$.}\label{fig:comp_cotas_er}
\end{figure}

Let us now consider large sparse Erd{\"o}s-R\'enyi graphs. In this case, there are no previous estimates available, neither of the degree-greedy algorithm or the independence number. We will thus only consider the results of our approximations, alongside simulations of the Glauber dynamics with large activation rate. The results are shown in Fig.\ \ref{fig:comp_cotas_er}. As expected, our approximation lays extremely close to the degree-greedy algorithm, which in turn is optimal for values of $\lambda$ smaller than 2.4 (note that the Glauber dynamics does not produce a larger independence number in these cases). The characterization provided by our approximation is an important contribution, as the exact value of the independence number was unknown for this family of random graphs. Finally, note that the static degree-aware algorithm (i.e.\ the relatively simpler to solve equation \eqref{eq:modulatedIS}) in this case produces a reasonable approximation to the independence number. 

\subsection{Geometric graphs}\label{sec:appgeom}

\begin{figure}
    \centering
    \includegraphics[width=0.6\textwidth]{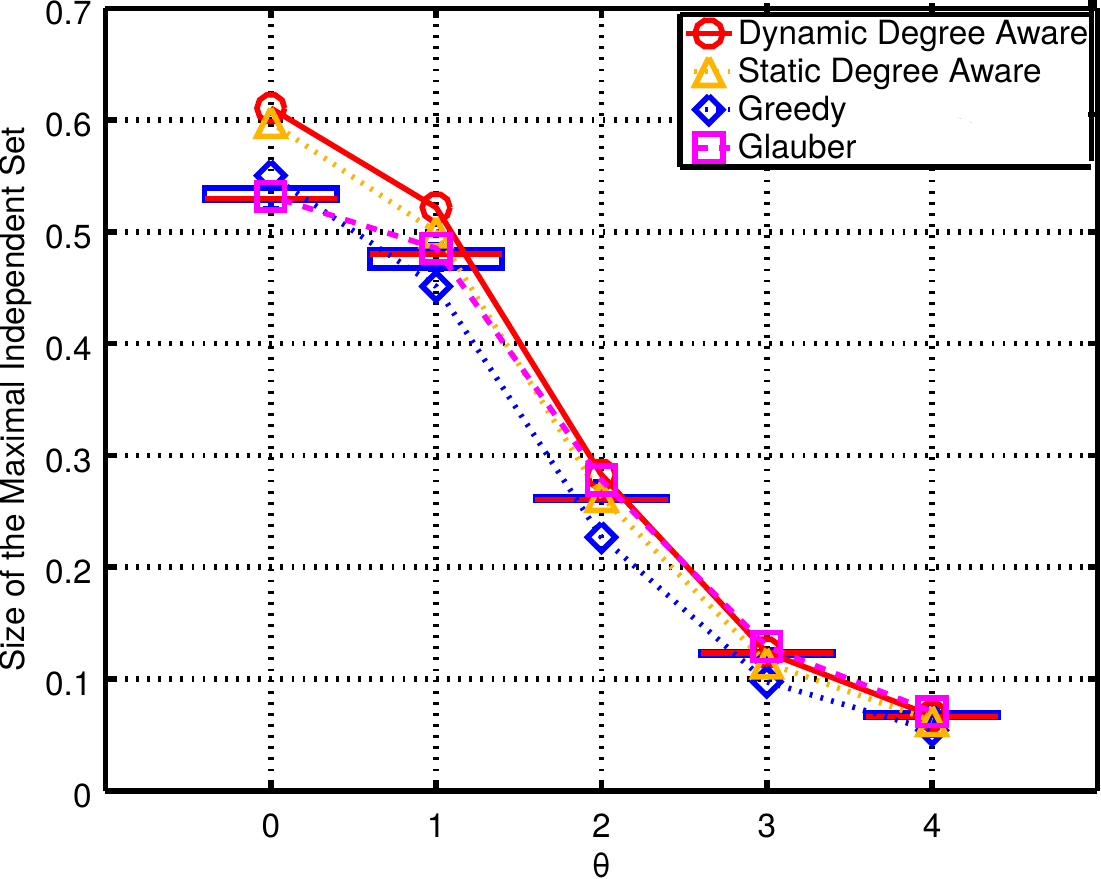}
    \caption{Size of the maximal independent set for different values of $\theta$ for the ppp case plus fading. }\label{fig:comp_poisson}
\end{figure}

We will now consider a case closely related to the Erd{\"o}s-R\'enyi model: graphs stemming from a Poisson Point process (ppp) on the plane. In these simulations, we generate a ppp on a circle of a proper size so that the mean number of nodes is fixed at $n=1000$. With the objective of modeling wireless networks, given two nodes $i$ and $j$, we will consider that they are connected if:
\begin{gather*}
    P(i,j) = X_{i,j}d(i,j)^{-a}>T,
\end{gather*}
where $d(i,j)$ is the distance between nodes $i$ and $j$, $X_{i,j}$ is a random variable that models fading (and which we will assume log-normally distributed with mean equal to $1$, variance $\theta$ and symmetric on $i,j$), $a=2$, and $T$ such that the mean number of neighbors when no fading is present is equal to 2. Results are shown in Fig.\ \ref{fig:comp_poisson}. 

Note that the case where $\theta=0$ corresponds to a variant of the so-called Mat\`ern Hard Core, and that the initial degree distribution is the same as the Erd{\"o}s-R\'enyi case (with $\lambda=2$ in this particular case). However, the fact that the Mat\`ern Hard Core process stems from a spatial graph generates correlations between the degrees of the nodes, which are not taken into account by the Configuration Model. This results in both our approximations (i.e.\ eqs.\ \eqref{eqsdificil} and \eqref{eq:modulatedIS}) overestimating the size of the independent set found by the degree-greedy algorithm. However, as the fading increases, these correlations become increasingly smaller, resulting in a better approximation by our method. 

It is interesting to note that the case where $\theta=1$ and $\theta=0$, the Glauber dynamics is not capable of improving the results of the degree-greedy algorithm. The resulting degree distribution for both cases is (empirically verified for $\theta=1$) a binomial one with mean less than 2.4. This seems to indicate that results such as the one proposed in Sec.\ \ref{sec:optimality} are valid for a wider type of random graphs, and not just those generated by the Configuration Model, a research line worth following in the future.

\section{Maximum Independent Sets and Wireless Networks}
\label{sec:maxywifi}

We will now focus our attention on the problem of computing (or approximating) the capacity of a large 802.11-based wireless network (or Wi-Fi as it is commercially known). This kind of networks implement a medium access protocol called CSMA/CA (Carrier Sense Multiple Access with Collision Avoidance) algorithm.

In a nutshell, CSMA/CA works in the following way. When a node wants to transmit a packet, it first draws a random backoff timer, which is decremented as long as none of its neighbors is transmitting. When the backoff finally reaches zero the packet is transmitted. We will assume that nodes are always trying to transmit packets, so this process is repeated indefinitely.

As shown in~\cite{laufer2016capacity}, when the mean backoff time is much smaller than the mean transmitting time (which, if the system is designed for efficiency, should be the case) then the number of transmitting nodes, after a long enough time, approaches that of the maximum independent set. Intuitively, when a node stops transmitting, its neighbors and that same node start competing for the medium. We thus have that any decrease by one on the active or transmitting nodes is rapidly followed by an increase in \emph{at least} one. In fact, if both the backoff timer and the transmission time are exponentially distributed, this process behaves as a continuous time version of the Glauber dynamics. 

The previous discussion means that the capacity of a wireless network may be approximated by our method. In this subsection, we will complement the analysis we performed in the previous subsection by simulating a CSMA/CA over two families of graphs and showing how the number of active (transmitting) nodes evolves over time and how fast it converges to the independence number. 

\subsection{Erd{\"o}s-R\'enyi graphs}

Let us focus on the Erd{\"o}s-R\'enyi case. As we showed in Sec.\ \ref{sec:optimality}, the degree-greedy algorithm produces a maximum independent set in the case of $\lambda \lesssim 2.4$. To further illustrate this result (cf.\ Fig.\ \ref{fig:comp_cotas_er}), we now present the evolution of CSMA/CA (where the backoff and transmission times are exponentially distributed with means equal to $1\times 10^{-6}$ and $1\times 10^{2}$ respectively) over an Erd{\"o}s-R\'enyi graphs with a mean number of neighbors equal to $\lambda$. 

\begin{figure}
    \centering
    \includegraphics[width=0.6\textwidth]{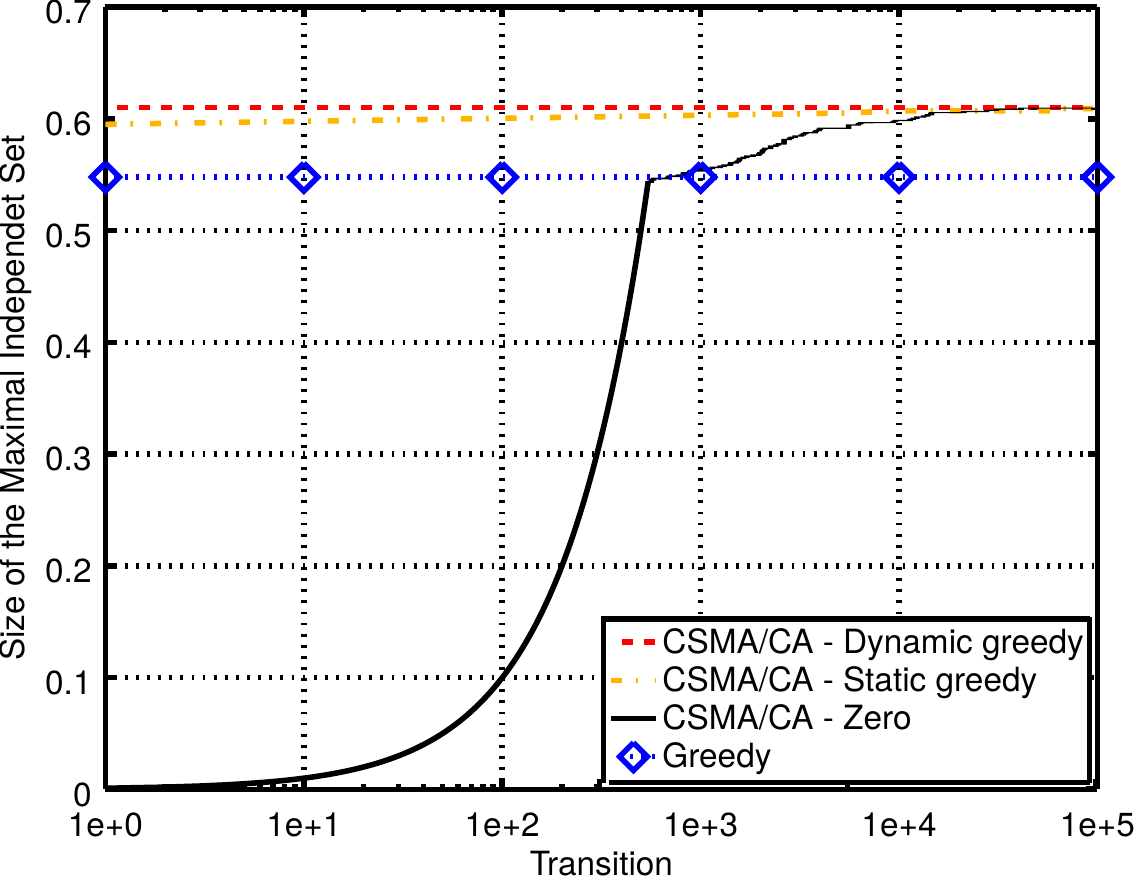}
    \caption{The number of transmitting nodes as the CSMA/CA algorithm proceeds (starting from three possible initial conditions) for an Erd{\"o}s-R\'enyi graph with $\lambda=2$. The abscissa in logarithmic scale.}\label{fig:comp_cotas_er_dcf}
\end{figure}

In particular, we ran the CSMA/CA algorithm starting from three different initial conditions: no nodes transmitting, a situation where nodes stemming from a static degree-aware exploration algorithm are transmitting, and one where those from a dynamic degree-aware algorithm are. The results corresponding to $\lambda=2$ are shown in Fig.\ \ref{fig:comp_cotas_er_dcf} (abscissa in logarithmic scale for easier visualization), along with the independence number corresponding to the greedy algorithm. 

Let us consider the case where all nodes start in the backoff stage. As discussed in Sec.\ \ref{sec:glauber}, the evolution of the CSMA/CA mechanism in this case may be separated into two stages. During the first stage, the CSMA/CA rapidly reaches a number of active nodes that is distributed as the greedy algorithm (since the CSMA/CA mechanism is designed to choose randomly among competing nodes). Once this maximal independent set is reached, CSMA/CA  is still able to increase the number of concurrent transmissions, but very slowly: it takes roughly the first 20.000 transitions to converge to a value close to the dynamic degree algorithm. 

It is also worth verifying how even after another further 80.000 transitions (and even when the initial configuration is modified to follow the static or dynamic degree aware algorithms) CSMA/CA is not able to increase the number of active nodes. This verifies that \eqref{eqsdificil} may be used to calculate the capacity of a wireless network when the condition of Proposition \ref{prop:condnearopt} is met. If this is not the case, it may yet prove very useful, as we further discuss in the following subsection. 

\subsection{Combining sequential algorithm and Glauber dynamics}

Let us then consider an Erd\"os-R\'enyi graph, but with a mean number of neighbors equal to $10$. In this case, as the conditions of Proposition \ref{prop:condnearopt} are not met, in order to estimate the independence number (or the capacity of the wireless network) one has to resort to algorithms such as Glauber (or equivalently CSMA/CA), which will eventually converge. However, based on the previous simulations and our results, to ``accelerate'' this estimation we propose to start a Glauber dynamics from a configuration given by one of the sequential algorithms we discussed (static or dynamic degree greedy). 

\begin{figure}
    \centering
    \includegraphics[width=0.6\textwidth]{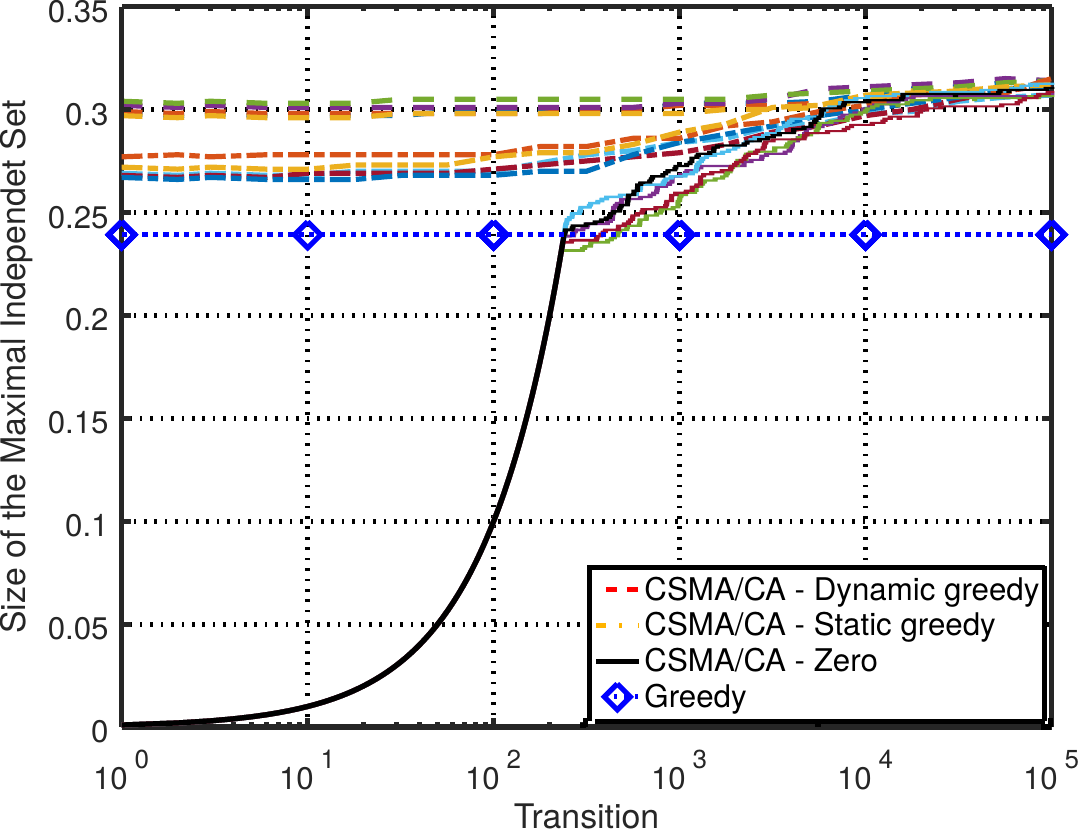}
    \caption{The number of transmitting nodes as the CSMA/CA algorithm proceeds (starting from three possible initial conditions) for an Erd{\"o}s-R\'enyi graph with $\lambda=10$. The abscissa in logarithmic scale.}\label{fig:comp_cotas_er_dcf_10vec}
\end{figure}

Simulation results are shown in Fig.\ \ref{fig:comp_cotas_er_dcf_10vec}, where we have included five executions of a CSMA/CA for each of the initial configurations (no nodes transmitting and both sequential algorithms). It is important to highlight that these initial configurations may be found in a number of steps equal to the number of nodes in the maximal independent set, and thus a fraction of the number of nodes in the graph ($n=1000$ in these simulations). Note that the trend observed in our previous simulations still holds: starting from no nodes transmitting, the number of transitions it takes CSMA/CA to reach an independent set of size similar to that of the dynamic degree greedy algorithm is (at least) an order of magnitude larger than the number of nodes in the graph. Moreover, even after 100.000 iterations, the difference among the initial configurations is not negligible, and the best results are obtained by starting the CSMA/CA algorithm from the configuration found by a dynamic degree greedy algorithm. 

\subsubsection{Geometric graphs} 

\begin{figure}
    \centering
    \includegraphics[width=0.6\textwidth]{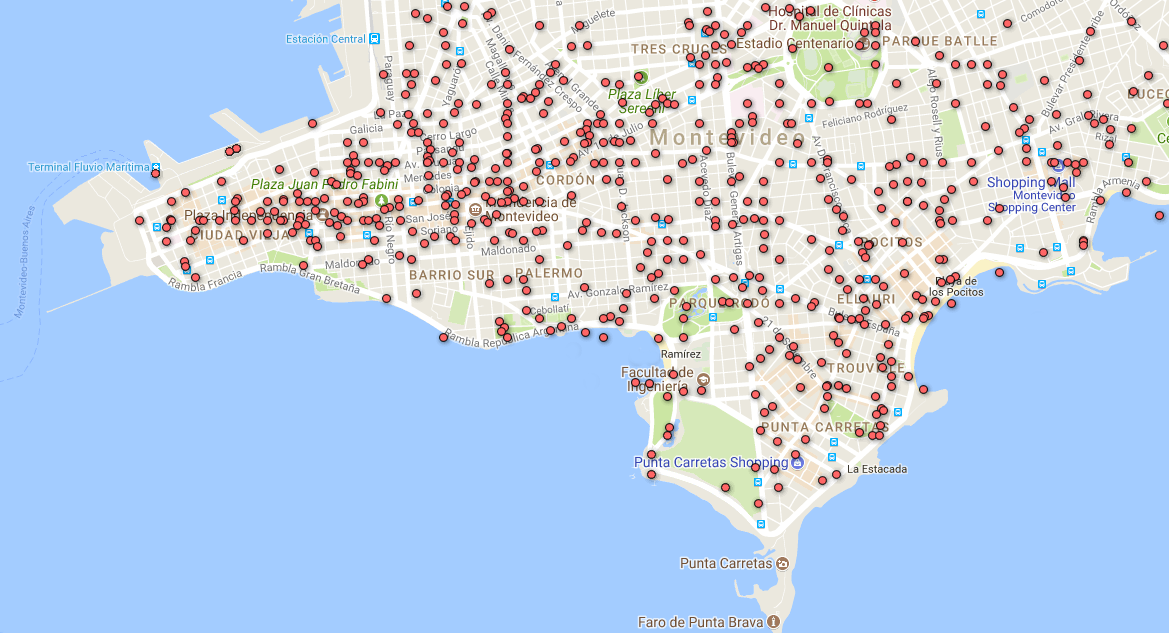}
    \caption{The nodes and their position.}\label{fig:mapa}
\end{figure}

Let us now discuss a more realistic case-scenario. Instead of arbitrarily positioning the nodes according to a ppp, we will consider their locations as provided by a public dataset (in this case OpenCellID~\cite{opencell}, an open and crowdsourced database of cell tower locations around the world). Figure \ref{fig:mapa} shows the area considered in this example, consisting of an urban region of around 20\,$km^2$ and including 579 nodes.

\begin{figure}
    \centering
    \includegraphics[width=0.6\textwidth]{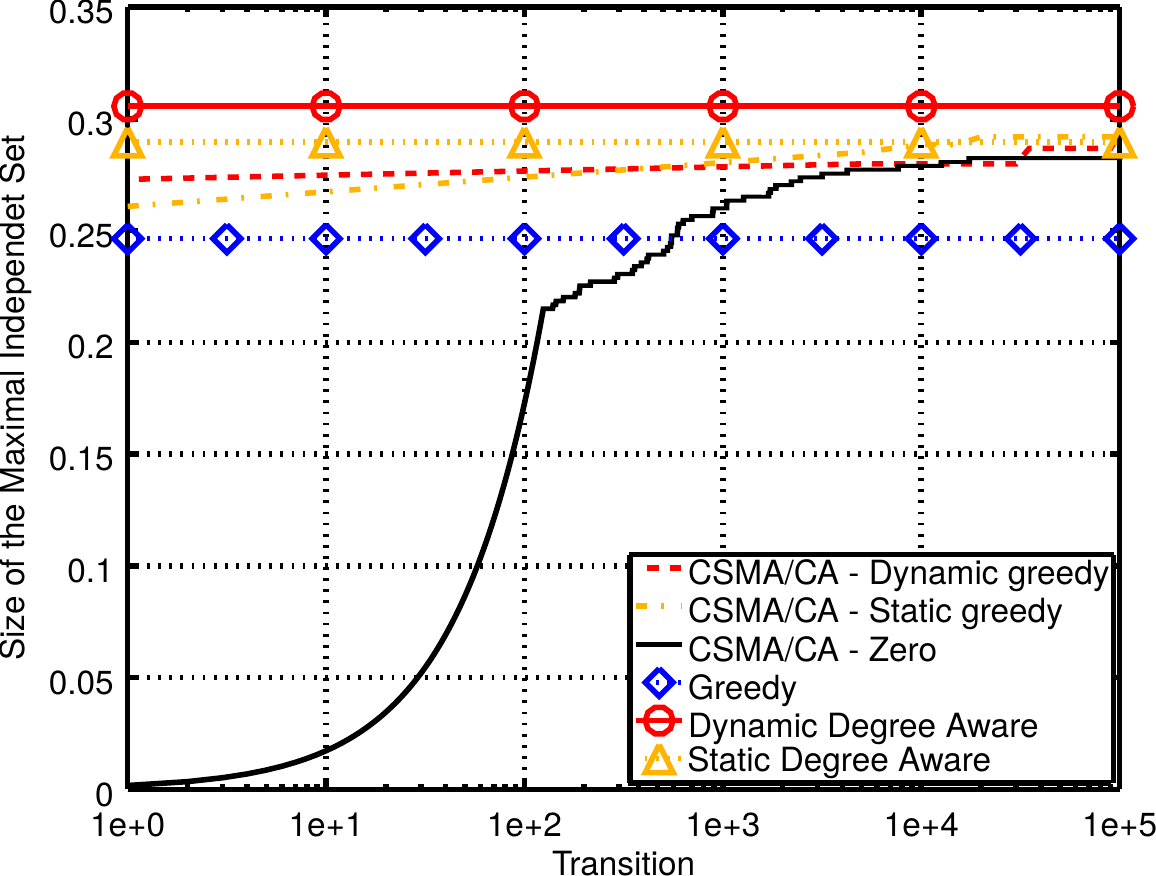}
    \caption{The number of transmitting nodes as the CSMA/CA algorithm proceeds (starting from three possible initial conditions) for the nodes shown in Fig.\ \ref{fig:mapa}. The abscissa in logarithmic scale. }\label{fig:results_montevideo}
\end{figure}

We consider the same propagation model as before, with log-normal fading with variance $\theta$ and an average number of neighbors equal to 2.4 when $\theta=0$ (equivalent to a range of 150\,m). Convergence results corresponding to $\theta=3$ are shown in Fig.\ \ref{fig:results_montevideo}, along with the three theoretical results.  Some observations are in order. First, if fading is non-negligible as in this case (and as in Fig.\ \ref{fig:comp_poisson}), our approximation based on configuration models is very precise. Further simulations, not reported here due to space limitations, indicate that this trend is true for other scenarios, such as sub-urban or dense-urban. This illustrates the flexibility of analyses based on Configuration Models in general, and of our approximation in particular. Secondly, and differently to the case of Fig.\ \ref{fig:comp_cotas_er_dcf}, the three initial conditions result in a non-negligible difference in the number of active nodes, even after 100.000 iterations. This observation further illustrates that Glauber dynamics can get stuck in local minima while starting them from advantageous configurations can lead to discover bigger independent sets much faster. 

\subsection{Fairness and independent set}\label{sec:fairness}

Until now, we have not discussed how \emph{unfair} is the dynamic that establishes the connections in the network; but in a greedy (and even more in a degree-greedy) kind of dynamics, vertices with a higher degree tend to be blocked faster than vertices with a smaller one. This is so because, if the degree of a vertex is large, it has a bigger chance of some of its neighbors activating before he does. Here we discuss how to compensate this by implementing a family of static degree-aware processes.

As a measure of the \emph{unfairness} we will use the \emph{Total Variation} $||\cdot||_{TV}$ between the degrees of the vertices of the independent set constructed by the dynamics $(q_i)_{i\geq0}$ and the initial degree distribution of the entire graph $(p_i)_{i\geq0}$, given by
\begin{equation}\label{eq:tv}
    ||p - q||_{TV} = \frac{1}{2} \sum_{i=1}^{\infty} |p_i - q_i|. 
\end{equation}

It is easy to show that when $||p - q||_{TV}=0$, the probability of every node of activating in some point of the evolution of the process is the same. On the other hand, higher values of Total Variation will correspond to situations in which nodes of some degrees connect with a higher probability than others; that is, with less fair situations.

Recall that to describe the limit of a static degree-aware exploration it is enough to integrate equation a single ODE. After doing so, the asymptotic proportion of vertices of the independent set found can be determined by \eqref{eq:modulatedIS} while the distribution $(q_i)_{i\geq0}$ can be obtained by \eqref{eq:modulateddegree}. 

\begin{figure}
    \centering
    \includegraphics[width=0.6\textwidth]{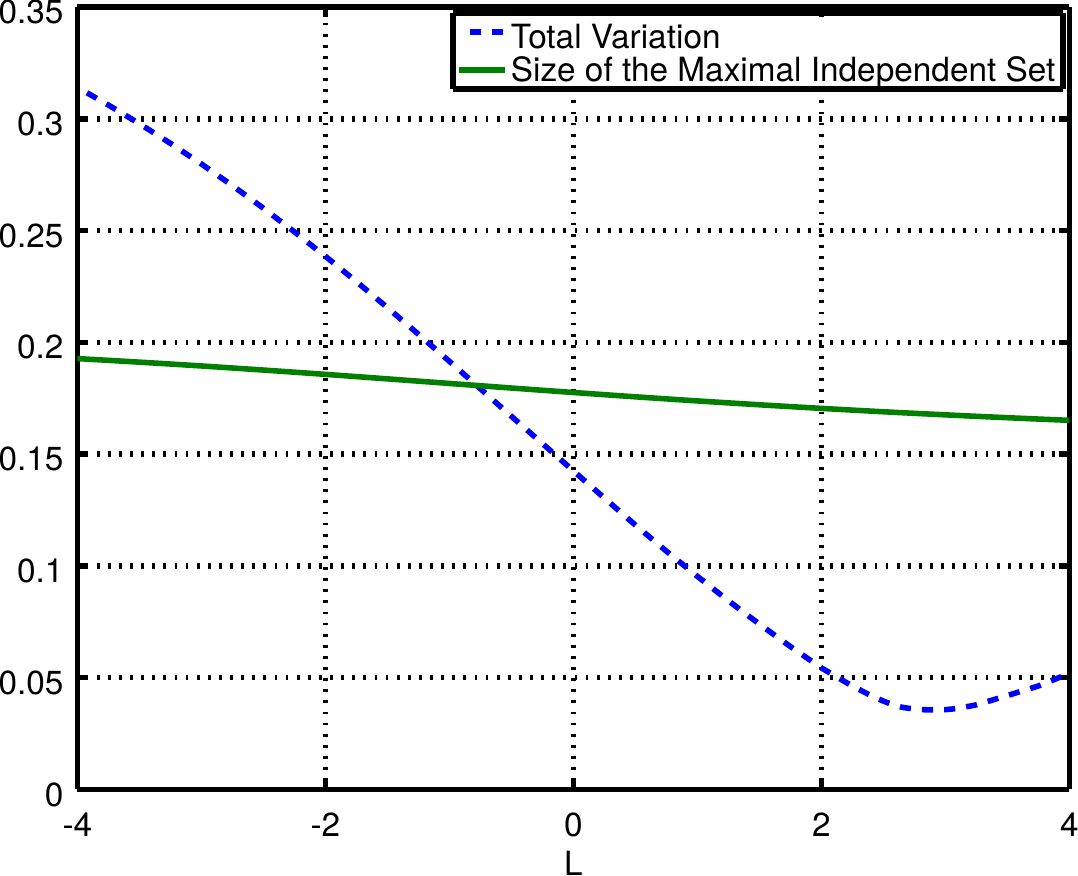}
    \caption{Unfairness and proportion of vertices in the independent set as a function of the power $L$ of the clock distribution for a Poisson degree distribution of mean 16.}\label{fig:simufairness}
\end{figure}

Using this relation, the unfairness of static degree-aware processes with different sets of clocks can be computed. As an example, in Fig.\ \ref{fig:simufairness}, we present the unfairness and proportion of vertices in the independent set for clocks given by $\lambda(i) = (i+1)^L$ with $L \in \mathbb{R}$ as a function of the power $L$ for a Poisson degree distribution of mean 16.

As can be seen in the figure, there is a trade-off between fairness and the independent set size. Negative powers obtain a larger independent set but result in higher unfairness than regular greedy dynamics (which corresponds to $L = 0$). On the other hand, higher fairness can be achieved by assigning clocks that depend on the degree as positive powers at the expense of a smaller independence set. Intuitively, this is so because with negative powers, vertices with lower degree tend to connect faster which results in less blocked vertices and ultimately in a larger independent set. On the other hand, positive powers make nodes with higher degree connect faster which tends to compensate the natural unfairness of the greedy type dynamics discussed at the beginning of the section. 

From the figure we can see that there is a minimum in unfairness in the interval $(2.5,3)$. This reduction in unfairness is achieved at the cost of a reduction of less than $10\%$ in the size of the independent set. This means that unfairness in transmissions can be considerably improved while not reducing much the number of transmitting nodes. A similar situation is observed for a range of Poisson distributions with mean between 4 and 30.

\section{Conclusions and future work}

We extended existing hydrodynamic limits over Configuration Model graphs to two variants of degree-aware algorithms: dynamic and static. The dynamic algorithm is capable of approximating to an arbitrary precision the degree-greedy algorithm. Thus, through a numerical evaluation of a system of differential equations, we estimate the jamming constant of degree-greedy exploration algorithms for much more general random graphs than previous results, in addition to characterizing the independence number of this family of graphs. 

The static dynamics not only provides closed form results for the size of the discovered independent set but also for the degree distribution of its nodes.  Although sub-optimal, this property allows us to analyze the trade off between fairness and the size of the corresponding independent set for different random graph distributions and to propose a parameter selection which ensures minimum unfairness along with a mild reduction on this size.  This result provides insights to the issue of fairness in wireless networks, also known as the ``starvation problem".

We have illustrated the usefulness of our results by estimating the capacity of a 802.11-based network, both for synthetic as well as real networks. Moreover, a combination of sequential and dynamical (Glauber) algorithms is proposed and analyzed in order to improve the time of convergence. It is interesting to note that previous results on the capacity of these networks have either focused on mean-field situations where all nodes ``see'' each other (for example in the seminal paper by Bianchi~\cite{bianchi2000performance}, to quote a single reference), or on small given graphs where calculating the independence number is feasible (such as in~\cite{soung2010back,laufer2016capacity}). Hence, random graph methods give a viable alternative to stochastic geometry, whose complexity can be prohibitive. 
For a more thorough discussion and comparison with stochastic geometry see~\cite{rattaro2017estimating}.

As future work, we plan to characterize theoretically the asymptotic optimality threshold in function of the degree distribution.

%
%
%




\begin{thebibliography}{}

\bibitem[{Bayati et al.(2010)}]{bayati2010combinatorial}
Bayati M, Gamarnik D, Tetali P (2010) Combinatorial approach to the interpolation method and scaling limits in sparse random graphs. \emph{Proceedings of the forty-second ACM symposium on Theory of computing}, 105--114.

\bibitem[{Bermolen et al.(2017a)}]{bermolen2017jamming}
Bermolen P, Jonckheere M, Moyal P (2017) The
jamming constant of uniform random graphs.
\emph{Stochastic Processes and their Applications} 127, 7, 2138–2178.

\bibitem[{Bermolen et al.(2017b)}]{Bermolen2017}
Bermolen P, Jonckheere M, Sanders J (2017) Scaling Limits and Generic Bounds for Exploration Processes. \emph{Journal of Statistical Physics} 169, 5, 989--1018.

\bibitem[{Bermolen et al.(2018)}]{PER}
Bermolen P, Larroca F, Jonckheere M, S\'aenz M (2018) Degree-Greedy Algorithms on Large Random Graphs. \emph{ACM SIGMETRICS Performance Evaluation Review} [\emph{To appear}]. 

\bibitem[{Bianchi(2000)}]{bianchi2000performance}
Bianchi G (2000) Performance analysis of the IEEE 802.11 distributed coordination function. \emph{IEEE Journal on Selected Areas in Communications} 18, 3, 535--547.

\bibitem[{Bollobas and (1976)}]{bollobas1976cliques}
Bollob{\'a}s B, Erd{\"o}s P (1976) Cliques in random graphs, \emph{Mathematical Proceedings of the Cambridge Philosophical Society} 80, 3, 419--427.

\bibitem[{Bollobas(1981)}]{bollobas1981independence}
Bollob\'as B (1981) The independence ratio of regular graphs, \emph{Proceedings of the American Mathematical Society}, 433--436.

\bibitem[{Bollobas(1998)}]{bollobas01}
Bollob\'as B. Random graphs, 215--252. Springer.

\bibitem[{Brightwell et al.(2017)}]{janson2017greedy}
Brightwell G, Janson S, Luczak M (2017) The
Greedy Independent Set in a Random Graph with
Given Degrees. \emph{Random Structures and Algorithms} 51, 4, 565–586.

\bibitem[{Ding et al.(2016)}]{ding2016maximum}
Ding J, Sly A, Sun N (2016) Maximum independent sets on random regular graphs. \emph{Acta Mathematica}, 217, 2, 263--340.

\bibitem[{Durrett(2007)}]{Durrett}
Durrett R (2007) Random graph dynamics. Cambridge university press Cambridge.

\bibitem[{Frieze(1990)}]{frieze1990independence}
Frieze A M (1990) On the independence number of random graphs. \emph{Discrete Mathematics} 81, 2, 171--175.

\bibitem[{Gamarnik and Goldberg(2010)}]{gamarnik2010randomized}
Gamarnik D, Goldberg D A (2010) Randomized greedy algorithms for independent sets and matchings in regular graphs: Exact results and finite girth corrections. \emph{Combinatorics, Probability and Computing} 19, 1, 61--85.

\bibitem[{Galvin and Tetali(2006)}]{galvin2006slow}
Galvin D, Tetali P (2006) Slow mixing of Glauber dynamics for the hard-core model on regular bipartite graphs. \emph{Random Structures \& Algorithms} 28, 4, 427--443.

\bibitem[{Halld{\'o}rsson and Radhakrishnan(1997)}]{halldorsson1997greed}
Halld{\'o}rsson M, Radhakrishnan J (1997) Greed is good: Approximating independent sets in sparse and bounded-degree graphs. \emph{Algorithmica} 18, 145--163.

\bibitem[{Jonckheere and S\'aenz(2018)}]{jonckheere2018asymptotic}
Jonckheere M, S\'aenz M (2018) Asymptotic opitmality of degree-greedy discovering of independent sets in Configuration Model graphs. arXiv preprint arXiv:1808.10358.

\bibitem[{Karp and Sipser(1981)}]{karp1981maximum}
Karp R M, Sipser M (1981) Maximum matching in sparse random graphs. \emph{Foundations of Computer Science, 1981. SFCS'81. 22nd Annual Symposium}, 364--375.

\bibitem[{Lauer and Wormald(2007)}]{lauer2007large}
Lauer J, Wormald N C (2007) Large independent sets in regular graphs of large girth. \emph{Journal of Combinatorial Theory, Series B} 97, 6, 999--1009.

\bibitem[{Laufer and Kleinrock(2016)}]{laufer2016capacity}
Laufer R, Kleinrock L (2016) The Capacity of
Wireless CSMA/CA Networks. \emph{IEEE/ACM Transactions on Networking} 24, 3 (June 2016), 1518--1532.

\bibitem[{Liew et al.(2010)}]{soung2010back}
Liew S C, Kai C H, Leung H C, Wong P (2010) Back-of-the-Envelope Computation of Throughput Distributions in CSMA Wireless Networks. \emph{IEEE Transactions on Mobile Computing} 9, 9 (Sept 2010), 1319--1331.

\bibitem[{McKay(1987)}]{mckay1987lnl}
McKay B D (1987) Independent sets in regular graphs of high girth. \emph{Ars Combinatoria} 23, 179--185.

\bibitem[{Molloy and Reed(1998)}]{molloy1998size}
Molloy M, Reed B (1998) The size of the giant component of a random graph with a given degree sequence. \emph{Combinatorics, probability and computing} 7, 3, 295--305.

\bibitem[{Opencell(2018)}]{opencell}
unwiredlabs (2018) OpenCellID. \url{http://opencellid.org/}

\bibitem[{Rattaro et al.(2017)}]{rattaro2017estimating}
Rattaro C, Larroca F, Bermolen P, Belzarena P (2017) Estimating the medium access probability in large cognitive radio networks. \emph{Ad Hoc Networks} 63, 1--13.

\bibitem[{van der Hofstad(2016)}]{remco2016}
van der Hofstad R (2016) Random Graphs and Complex Networks: Volume 1. Cambridge University Press.

\bibitem[{Vigoda(2001)}]{vigoda2001note}
Vigoda E (2001) A note on the Glauber dynamics for
sampling independent sets. \emph{The electronic journal of combinatorics} 8, 1.

\bibitem[{Wormald(1995)}]{wormald1995differential}
Wormald N C (1995) Differential equations for random processes and random graphs. \emph{The annals of applied probability}, 1217--1235.

\bibitem[{Wormald(1999)}]{wormaldCM}
Wormald N C (1999) Models of random regular graphs. \emph{London Mathematical Society Lecture Note Series}, 239--298. Cambridge University Press.


\end{thebibliography}



\bibliographystyle{nonumber}

\end{document}